\documentclass[oneside,openany,11pt]{article}

\title{Gromov-Witten invariants of Fano hypersurfaces, revisited}
\author{Hironori Sakai}
\date{}

\usepackage{amsmath,amssymb} 

\newtheorem{df}{Definition}[section]
\newtheorem{prop}[df]{Proposition}
\newtheorem{thm}[df]{Theorem}
\newtheorem{lem}[df]{Lemma}
\newtheorem{cor}[df]{Corollary}

\usepackage[slantedGreek,noBBpl]{mathpazo}

\def\Z{{\mathbb Z}}

\def\C{{\mathbb C}}

\def\D{{\mathcal D}}

\def\L{{\mathcal L}}
\def\M{{\mathcal M}}

\DeclareMathOperator{\End}{End}

\DeclareMathOperator{\qh}{QH}
\DeclareMathOperator{\gw}{GW}
\DeclareMathOperator{\diag}{diag}

\newcommand{\proof}[0]{\medskip\noindent \textit{Proof.}\quad}
\newcommand{\proofof}[1]{\medskip\noindent \textit{Proof of #1.}\ }
\newcommand{\qed}[0]{\hfill $\square$\bigskip}

\newcommand{\rd}[0]{\partial}

\newcommand{\pair}[1]{\langle#1\rangle}

\newcommand{\cp}[1]{\C P^{#1}}

\newcommand{\hs}[2]{H^{\sharp}(M_{#1}^{#2})}
\newcommand{\qs}[2]{\mathrm{QH}^{\sharp}(M_{#1}^{#2})}

\newcommand{\PF}[0]{\Omega^h_\mathrm{PF}}
\newcommand{\dub}[0]{\nabla^h_\mathrm{D}}
\newcommand{\DUB}[0]{\Omega^h_\mathrm{D}}
\newcommand{\JIN}[0]{\Omega^h_\mathrm{J}}

\newcommand{\agg}[0]{\mathcal{G}_\mathrm{AD}}
\newcommand{\afc}[0]{\mathcal{A}_\mathrm{AD}}

\begin{document}
\maketitle

\begin{abstract}
 The goal of this paper is to give an efficient computation of the
 3-point Gromov-Witten invariants of Fano hypersurfaces,
 starting from the Picard-Fuchs equation. This simplifies and
 to some extent explains the original computations of Jinzenji.
 The method involves solving a gauge-theoretic differential equation,
 and our main result is that this equation has a unique solution.
\end{abstract}

\section{Introduction}

Gromov-Witten invariants compute
``numbers of pseudo-holomorphic curves'' in a symplectic manifold.
They are rigorously defined as integrals on moduli spaces of stable maps.
Therefore it is difficult to calculate Gromov-Witten invariants
directly from the definition.

An alternative method of computation comes from mirror symmetry.
Although the mirror symmetry phenomenon is not yet completely understood,
it suggests that Gromov-Witten invariants can be computed in terms of
coefficients of power series solutions of certain differential equations.
The most well known example is the quintic hypersurface in $\cp{4}$;
this is a Calabi-Yau 3-fold. Fano hypersurfaces are more elementary
from the point of view of Gromov-Witten invariants, and it was 
established by Givental that the Gromov-Witten invariants in this 
case are determined by the ``Picard-Fuchs equation''.
The Picard-Fuchs equation of the quintic hypersurfaces in $\cp{4}$ is
$\bigl( \rd^4 - 5e^t(5\rd+4)(5\rd+3)(5\rd+2)(5\rd+1) \bigr) \psi(t)=0$.

A hypersurface $M^k_N$ of degree $k$ in $\cp{N-1}$ is
Fano if and only if $N > k$, and the Picard-Fuchs equation is
\[
 \Bigl( \rd^{N-1} - ke^t(k\rd + (k-1)) \dots (k\rd +2)(k\rd + 1) \Bigr)
 \psi(t)=0.
\]
Before Givental's work, partial results on the quantum cohomology of 
Fano hypersurfaces had been obtained by Collino-Jinzenji
(\cite{collino99}) and Beauville (\cite{beauville95:_quant}).
Subsequently, Jinzenji (\cite{jinzenji02:_gauss_manin}) observed
that a simple Ansatz leads to the correct Gromov-Witten invariants
and he obtained complicated but explicit formulae from this Ansatz.

The aim of this paper is to calculate the $3$-point Gromov-Witten
invariants of a Fano hypersurface by using the method of
\cite{amarzaya:_gromov_witten_d,guest05:_quant_d}.
In this method, the flat connection associated to the $\D$-module
$\D/(\mathrm{PF})$ is ``normalized'' by applying the Birkhoff 
factorization.  We shall show (as a consequence of Givental's work) that
this method produces the correct $3$-point Gromov-Witten invariants.

The algorithm for the computation of $3$-point Gromov-Witten invariants
from the quantum differential equations was introduced in
\cite{guest05:_quant_d}, and applied to flag manifolds in
\cite{amarzaya:_gromov_witten_d}, and our treatment of hypersurfaces is
broadly similar. However, there are some special features in this case,
which makes a separate discussion worthwhile.  First, the differential 
equations in this case is are o.d.e., rather than p.d.e.,
hence the integrability condition plays no role. Second, the
o.d.e. which appears in the Birkhoff factorization can be integrated
very explicitly, and this leads to purely algebraic formulae (whereas
the algorithm in \cite{amarzaya:_gromov_witten_d} required the solution
of large systems of p.d.e.).

Computationally, our method is similar to Jinzenji's method,
but considerably simpler.
In \S $2$ we review the definitions of the Gromov-Witten invariants,
the quantum cohomology ring, and the Dubrovin connection.
In \S $3$ we discuss the quantum cohomology of Fano hypersurfaces.
The quantum differential system and an example of Jinzenji's results
are discussed in \S $4$. 
In \S $5$, we explain the loop group method and we compute a flat
connection from a $\D$-module which is related to the quantum
differential system.
In \S $6$, we discuss relation between families of connection $1$-forms
and $\D$-modules. The ``adapted'' gauge group is the most important
object.
In \S $7$, we explain Jinzenji's results from our viewpoint and
prove that our results agree with Jinzenji's results.
We also prove that our results produce the Gromov-Witten invariants.

\textbf{Acknowledgments}\quad
I would like to express my heartfelt thanks to Martin Guest, my adviser,
for his guidance and insight.

\section{Gromov-Witten invariants}

For any smooth projective manifold $M$ and any homology class
$A \in H_2(M;\Z)$, we can define the ($3$-pointed genus zero)
Gromov-Witten invariant $\gw_A : H^*(M;\C)^{\otimes 3} \rightarrow \C$.
Let $X_i$ be a smooth submanifold of $M$ and $a_i \in H^*(M;\Z)$ the
Poincar\'{e} dual of the fundamental class $[X_i]$, where $i = 0, 1, \infty$.
Then the invariant $\gw_A(a_0,a_1,a_\infty)$ is roughly
the number of holomorphic maps $f:\cp{1} \rightarrow M$ such that
\[
 f_*[\cp{1}] = A \quad\mbox{and}\quad
 f(i) \in X_i \ \ \mbox{for}\ \ i = 0,1,\infty \in \cp{1}.
\]
For a rigorous definition, we must construct a certain moduli space.
The Gromov-Witten invariants satisfy properties which are called
Kontsevich-Manin axioms (see \cite{d.a.cox99:_mirror_symmet_algeb_geomet}
or \cite{mcduff04:_j} for details). In this paper we use the following
properties.
\begin{description}
 \item[Linearity Axiom.] $\gw_A$ is linear in each variable.
 \item[Effectivity Axiom.] Let $\omega$ be the symplectic form of $M$.
	    Then the Gromov-Witten invariant
	    $\gw_A$ vanishes whenever $\int_A \omega < 0$.
 \item[Degree Axiom.] Let $x$, $y$ and $z$ be
	    homogeneous cohomology classes.
	    Then $\gw_A(x,y,z) = 0$ unless
	    \[
	     \deg x + \deg y + \deg z = 2\dim_\C M + 2c_1(A).
	    \]
	    Here $c_1$ is the first Chern class of $M$.
 \item[Equivariance Axiom.] For any homogeneous cohomology classes
	    $x$, $y$ and $z$ of $M$,
	    \begin{align*}
	     \gw_A(x,y,z) &= (-1)^{\deg x + \deg y} \gw_A(y,x,z) \\
	     &= (-1)^{\deg y + \deg z} \gw_A(x,z,y).
	    \end{align*}
 \item[Point Mapping Axiom.] For any cohomology classes
	    $x$, $y$ and $z$,
	    \[
	    \gw_0(x,y,z)
	    = \int_M x \cup y \cup z.
	    \]
\end{description}

For any $t \in H^2(M;\C)$, we define the (small) quantum product
$\circ_t$ of $H^*(M;\C)$. Let $x, y \in H^*(M;\C)$ and $(\ ,\ )$ the
intersection pairing:
\[
 (x,y) := \int_M x \cup y.
\]
Since the pairing is non-degenerate, there is a unique cohomology class
$x \circ_t y \in H^*(M;\C)$ such that for any $z \in H^*(M;\C)$
\[
 (x \circ_t y,z) = \sum_{A \in H_2(M;\Z)}\gw_A (x,y,z)q^A,
\]
where $q=\exp(t(A))$. The above sum is finite if the variety $M$ is
Fano. The quantum product $\circ_t$ makes $H^*(M;\C)$ a ring
$\qh^*(M)$ with identity. The ring $\qh^*(M)$ is called the quantum
cohomology ring. $\qh^*(M)$ is supercommutative, namely, 
if $x, y \in H^*(M;\C)$ are homogeneous, then
\[
 x \circ_t y = (-1)^{\deg x \deg y} y \circ_t x.
\]

Let $b_0,\dots,b_s$ be a basis of $H^*(M;\C)$ and $b^0,\dots,b^s$ its
dual basis with respect to the intersection pairing. We can write an
explicit formula for the quantum product:
\[
 x \circ_t y = \sum_i(x \circ_t y,b_i)b^i
 = \sum_i \sum_{A \in H_2(M;\Z)} \gw_A (x,y,b_i)q^A b^i.
\]

Using the quantum product, we define a family of connection $\dub$
parameterized by $h \in \C^{\times}$ on the trivial bundle
$H^*(M;\C) \times H^2(M;\C) \rightarrow H^2(M;\C)$. The connection
$1$-form $\DUB = \dfrac{1}{h}A \in \Omega^1(H^2(M;\C),\End(H^*(M;\C)))$
is defined by
\[
 A_t : H^2(M;\C) \rightarrow \End(H^*(M;\C));\ A_t(x)y = x \circ_t y.
\]
The connection $\dub = d - \DUB$ is called the Dubrovin
connection. The connection $\dub$ is flat (see \cite{mcduff04:_j}).

\section{Fano hypersurfaces of complex projective spaces}

Let $N \geq 5$ be an integer. $M_N^k$ denotes a hypersurface of degree
$k$ in the projective space 
$\cp{N-1}$ and $\iota : M_N^k \rightarrow \cp{N-1}$ 
denotes the inclusion map. 
The hypersurface $M_N^k$ is connected and simply connected.
Let $p \in H^2(\cp{N-1};\Z)$ be the
hyperplane class and $b = \iota^* p$.
The $i$-times cup product of $b$ is denoted by $b_i$,
whereas the $i$-times quantum product is denoted $b^i$.
According to the hard Lefschetz and Lefschetz hyperplane theorems,
the cohomology ring $H^*(M_N^k;\C)$ consists of two parts:
\begin{description}
 \item[K\"{a}hler part $\hs{N}{k}$.] This is the space generated by the
	    K\"{a}hler form of
	    $M_N^k$. $\hs{N}{k} = \bigoplus_{i=0}^{N-2}\C b_i$.
 \item[Primitive part $P(M_N^k)$.] This is the orthogonal complement of
	    the K\"{a}hler part $\hs{N}{k}$ with respect to the
	    intersection pairing. $P(M_N^k)$ is a subspace of
	    $H^{N-2}(M_N^k;\C)$.
\end{description}

The first Chern class $c_1(M_N^k)$ is $(N-k)b$. We will only consider
the case $N>k$, i.e. the case where $M_N^k$ is Fano.

Now we are going to study the quantum cohomology ring of $M_N^k$.
\begin{thm}
 Let $\tilde{b}$ be the homology class defined as follows:
 \[
 \tilde{b} =
 \begin{cases} b & (N-k \geq 2), \\ b+k!q & (N-k=1). \end{cases}
 \]
 Then for any $a, b \in P(M_N^k)$ we have
 \begin{enumerate}
  \item $\tilde{b}^{N-1} = k^k q \tilde{b}^{k-1}$,
  \item $\tilde{b} \circ a = 0$, and
  \item $a \circ \tilde{b} = \dfrac{1}{k}(a,\tilde{b})(\tilde{b}^{N-2}-k^k q \tilde{b}^{k-2})$.
 \end{enumerate}
\end{thm}
(See \cite{collino99} for the proof.)
The second identity implies that $\hs{N}{k}$ is closed
with respect to the quantum product.
We denote this subring by $\qs{N}{k}$.
We concentrate on the quantum product of the subring $\qs{N}{k}$. 

Note that $b_0=1, b_1, \dots ,b_{N-2}$ form a basis of $\qs{N}{k}$ over
$\C$. Its dual basis with respect to the intersection pairing is
$\dfrac{1}{k}b_{N-2}, \dots, \dfrac{1}{k}b_0$, because
\[
 (b_i,b_j) = (b_{i+j},[M_N^k]) = ((\iota^* p)^{i+j},[M_N^k])
 = (p^{i+j},k[\cp{N-2}]) = k\delta^{N-2}_{i+j}.
\]
Thus for $0 \leq m \leq N-2$
\[
 b \circ b_{N-2-m}
 = \sum_j \sum_d \gw_{dA}(b,b_{N-2-m},b_j)q^{dA}\dfrac{1}{k}b_{N-2-j},
\]
where $A$ is a generator of $H_2(M_N^k;\Z)$.
$\gw_{dA}(b,b_{N-2-m},b_j)$ vanishes unless
\[
 \deg b + \deg b_{N-2-m} + \deg b_j
 = 2\dim_{\C}M_N^k + \pair{c_1(M_N^k),dA},
\]
namely $j=m-1+d(N-k)$. For a positive integer $d$, we define
$L^d_m$ by
\[
 L_m^d = \dfrac{1}{k}\gw_{dA}(b,b_{N-2-m},b_{m-1+d(N-k)}).
\]
Then we can write
\[
 b \circ b_{N-2-m} =
 b_{N-1-m} + \sum_{d \geq 1} L_m^dq^{d}b_{N-1-m-d(N-k)}
\]
where we abbreviate $q^{dA}$ to $q^{d}$. Note that $L_m^d$
vanishes unless
\[
 0 \leq m-1 + d(N-k) \leq N-2,\ 0 \leq N-1-m - d(N-k) \leq N-2
\]
i.e., $0 \leq m \leq (N-1) - (N-k)d$.

\section{The quantum differential system}

We consider the quantum differential system for $M_N^k$
which is defined as
\begin{align*}
 &\dfrac{\rd \psi_{N-2-m}}{\rd t} = \psi_{N-1-m}(t)
 + \sum_{d \geq 1} L_m^d e^{dt} \psi_{N-1-m-d(N-k)}(t),\\
 &\dfrac{\partial \psi_{N-2}}{\partial t} = 
 \sum_{d \geq 1} L_0^d e^{dt} \psi_{N-1-d(N-k)}(t),
\end{align*}
(where $m = 1, \dots, N-2$).
The following important fact was proved by Givental
in \cite{givental96:_equiv_gromov_witten}.
\begin{thm}[Givental]
 \label{thm:Givental}
 The Gauss-Manin system can be reduced to the Picard-Fuchs equation:
\[
 \Bigl( \rd^{N-1} - ke^t(k\rd + (k-1)) \dots (k\rd +2)(k\rd + 1) \Bigr)
 \psi_0(t)=0,
\]
 where $\rd$ means $\dfrac{\rd}{\rd t}$.
\end{thm}

Jinzenji proposed a method to complete the Gromov-Witten invariants of
Fano hypersurfaces from the Picard-Fuchs equation
in \cite{jinzenji02:_gauss_manin}.
We explain briefly Jinzenji's method for the Fano hypersurface $M_5^3$.

The Gauss-Manin system for $M_5^3$ is
\begin{align*}
 \frac{\rd \psi_0}{\rd t} &= \psi_1(t),\\
 \frac{\rd \psi_1}{\rd t} &= \psi_2(t) + L_2^1e^t \psi_0(t), \\
 \frac{\rd \psi_2}{\rd t} &= \psi_3(t) + L_1^1e^t \psi_1(t), \\
 \frac{\rd \psi_3}{\rd t} &= \hspace{26.7pt} + \hspace{2.75pt}
 L_0^1e^t \psi_2(t) + L_2^2e^{2t}\psi_0(t).\\
\end{align*}
Reducing this system, we have
\[
\Bigl( \rd^4 - e^t \left( (L_2^1 + L_1^1 + L_0^1)\rd^2 + (2L_2^1+L_1^1)\rd + L_2^1 \right) - e^{2t}(L_2^2-L_0^1L_2^1) \Bigr)\psi_0 = 0. 
\]
Givental's theorem implies that this differential equation becomes
\[
 \Bigl( \rd^4 - 3e^t(3\rd+2)(3\rd+1) \Bigr)\psi_0 = 0.
\]
Thus we conclude
\[
 L_0^1 = 6, \quad L_1^1 = 15, \quad L_2^1 = 6, \quad L_0^2 = 36.
\]

\section{Birkhoff factorization}\label{section:Birkhoff}

We modify the quantum differential system with the parameter $h$ as follows:
\begin{align*}
 &h\dfrac{\rd \psi_{N-2-m}}{\rd t} = \psi_{N-1-m}(t)
 + \sum_{d \geq 1} L_m^d e^{dt} \psi_{N-1-m-d(N-k)}(t),\\
 &h\dfrac{\partial \psi_{N-2}}{\partial t} = 
 \sum_{d \geq 1} L_0^d e^{dt} \psi_{N-1-d(N-k)}(t),
\end{align*}
(where $m = 1, \dots, N-2$).
We shall use the following modification of theorem \ref{thm:Givental}:
\begin{prop}
 \label{prop:reduction}
 The quantum differential system with the parameter $h$ can be reduced to
 the Picard-Fuchs equation with parameter $h$:
 \[
 \Bigl( (h\rd)^{N-1} - kqh^{k-1}(k\rd + (k-1)) \dots (k\rd +2)(k\rd + 1) \Bigr)
 \psi_0(t)=0,
 \]
 where $q=e^t$.
\end{prop}
We will prove the proposition in section \ref{section:Jinzenji}.

In this section and the next section, we will show that the Picard-Fuchs
equation with $h$ gives the Gromov-Witten invariants by using the
Birkhoff factorization, as in
\cite{guest05:_quant_d,amarzaya:_gromov_witten_d}.

Let $\Lambda = \C[q]$ be the polynomial ring generated by $q$ and
$\D$ be the module generated by $h\rd$ over $\Lambda(h)$.
First of all we consider the $\D$-module $\M^h = \D/(P^{N,k})$, where 
$(P^{N,k})$ is the left ideal generated by the operator
\[
 P^{N,k} = (h\rd)^{N-1} - kqh^{k-1}(k\rd+(k-1)) \dots (k\rd+2)(k\rd+1).
\]

Second, we introduce a family of (flat) connection $1$-forms
$\PF = \dfrac{1}{h}R^h(q)dt$. Put $P_0=1$, $P_1 = h\rd$, \dots,
$P_{N-2} = (h\rd)^{N-2}$. Then the equivalence classes
$[P_0], \dots, [P_{N-2}]$ form a $\Lambda(h)$-basis of the $\D$-module
$\M^h$. We define $R^h(q)$ by
\[
 h\rd([P_0],\dots,[P_{N-2}])
 = ([P_0],\dots,[P_{N-2}])R^h(q).
\]
Then $\PF = \dfrac{1}{h}R^h(q)dt$ is of the form
\[
 \PF = \frac{1}{h}\omega + \theta_0 + h\theta_1 + \dots + h^{k-2}\theta_{k-2},
\]
where $\omega, \theta_0, \dots, \theta_p$ are matrix-valued $1$-forms
independent of $h$.

Finally we obtain a connection from $\PF$ by using the Birkhoff
decomposition which is a candidate for the Dubrovin connection.
We consider $h$ as a parameter in $S^1 \subset \C$. Since $\PF$ is flat,
there is a map $L$ from an open subset $V$ of $\C$ to the loop group  
$\Lambda GL(\C^{N-1})$ such that $\PF = L^{-1}dL$. The loop group
$\Lambda GL(\C^{N-1})$ is the group of all smooth map from $S^1$ to
$GL(\C^{N-1})$. 

Let $L=L_-L_+$ be the Birkhoff decomposition of $L$, where $L_-$ extends
holomorphically to $1 < |h| \leq \infty$ and $L_+$ to $|h| < 1$, and
$L_-|_{h=0}=I$. In other words, $L_-$ and $L_+$ have expansions in $h$ as
follows:
\[
 L_- = I + \frac{1}{h}A_1 + \frac{1}{h}A_2 + \dots,\quad
 L_+ = Q_0(I + h Q_1 + h^2 Q_2 + \dots).
\]
Note that the Birkhoff decomposition exists if and only if $L$ takes
values in the big cell of the loop group. For any $q \in V$ there
exists $\gamma \in GL(\C^{N-1})$ such that $\gamma L$ belongs to
the big cell around $q$. Thus we may assume that there is a Birkhoff
decomposition of $L$. 

If we expand $\hat\Omega^h = (L_-)^{-1}dL_-$ as a series in $h$,
then only negative powers of $h$ appear. On the other hand, we have
\begin{align*}
 (L_-)^{-1}dL_-
 &= (LL_+^{-1})^{-1}d(LL_+^{-1}) \\
 &= (L_+L^{-1})\left((dL)L_+^{-1} + Ld(L_+^{-1})\right) \\
 &= L_+ (L^{-1}dL)L_+^{-1} + L_+d(L_+^{-1}) \\
 &= L_+ \left(\frac{1}{h}\omega + \theta_0 + \dots + h^p\theta_p \right)L_+
 + L_+d(L_+^{-1}). 
\end{align*}
Since the negative powers of $h$ disappear except for
$\dfrac{1}{h}Q_0 \omega Q_0^{-1}$ in the above expression,
we conclude that $\hat\Omega^h = \dfrac{1}{h}Q_0 \omega Q_0^{-1}$.
In the next section we will see that $\hat\nabla^h = d + \hat\Omega^h$
agrees with the restricted Dubrovin connection $\dub$, where the restricted
Dubrovin connection is the restriction of the Dubrovin connection 
to the trivial bundle $\hs{N}{k} \times H^2(M_N^k;\C) \rightarrow \hs{N}{k}$.
The restriction is well-defined because the quantum product on $\hs{N}{k}$
is closed.

It is difficult to execute the Birkhoff decomposition in general.
Note that we need only $L_+$ to work out $\hat\Omega^h$.
Since the non-negative powers of $h$ in $(L_-)^{-1}dL_-$ disappear,
we can obtain differential equations for $L_+$:
\begin{prop}[\cite{amarzaya:_gromov_witten_d}]\label{prop:Lplus}
 $L_+ = Q_0(I + hQ_1 + h^2Q_2 + \dots)$ satisfies the following
 differential equations:
 \begin{align*}
  (\L_0) &\quad dQ_0 = Q_0(\theta_0 + [Q_1,\omega]),\\
  (\L_1) &\quad dQ_1 = \theta_1 + [Q_1,\theta_0] + [Q_2,\omega]
  - [Q_1,\omega]Q_1,\ and\\
  (\L_i) &\quad dQ_i = \theta_i + Q_1\theta_{i-1} + \dots + Q_{i-1}\theta_1
  + [Q_i,\theta_0] + [Q_{i+1},\omega] - [Q_1,\omega]Q_i
 \end{align*}
 for $i \geq 2$.
\end{prop}
Here, $\L_i$ denotes the equation for $Q_i$.

To calculate $L_+$ we introduce some notation.
Let $E_{i,j}$ be the $N-1 \times N-1$ matrix
with $(i,j)$-component $1$ and all other components zero.

For an integer $n$ with $|n| \leq N-1$, we define an $N-1 \times N-1$ matrix
$\diag_n(a_1,\dots,a_{N-2-|n|})$ by
\[
 \diag_n(a_1,\dots,a_{N-2-|n|}) =
 \begin{cases}
  \sum_{i=1}^{N-2-n} a_i E_{i,n+i} \quad (n \geq 0), \\
  \sum_{i=1}^{N-2+n} a_i E_{i-n,i} \quad (n < 0).
 \end{cases}
\]
We call this an $n$-diagonal matrix. The identity
$E_{i,j}E_{\alpha,\beta} = \delta_{j,\alpha}E_{i,\beta}$ implies
that the product of an $n$-diagonal matrix and an $m$-diagonal matrix is 
an $(n+m)$-diagonal matrix.

The matrix $\diag_{-1}(1,\dots,1)$ is denoted by $I_{-1}$. For a matrix 
$A=(a_{i,j})$ and non-negative integer $n$, we call the matrix
$\diag_n(a_{1,1+n},a_{2,2+n},\dots,a_{N-2-n,N-2})$
the $n$-diagonal component of $A$.

Furthermore we define non-negative integers $\lambda_i^k$ as in the
previous section by
\[
 k \prod_{j=1}^{k-1}(kX+j)
 = \lambda_{k-1}^{k}X^{k-1} + \lambda_{k-2}^{k}X^{k-2}
 + \dots + \lambda_1^k X + \lambda_0^k.
\]
The Picard-Fuchs operator of $M_N^k$ is described in terms of
$\lambda_i^k$ as follows:
\[
  P^{N,k} = (h\rd)^{N-1}
 - q\Bigl(
 \lambda_{k-1}^k(h\rd)^{k-1} + \lambda_{k-2}^kh(h\rd)^{k-2}
 + \dots + \lambda_0^kh^{k-1}
 \Bigr) .
\]
If there is no danger of confusion, we omit the upper suffix of $\lambda_i^k$.

Recall that $P_i = (h\rd)^i$. We have:
$[\rd P_0] = \dfrac{1}{h}[P_1], \dots, [\rd P_{N-3}] = \dfrac{1}{h}[P_{N-2}]$
and
\[
 [\rd P_{N-2}] 
 = \dfrac{1}{h}[(h\rd)^{N-1}]
 = \dfrac{1}{h}\lambda_{k-1}q[P_{k-1}]
 + \lambda_{k-2}q[P_{k-2}] + \dots + h^{k-2} \lambda_0q[P_0].
\]
Thus
\[
 \PF =
 \dfrac{1}{h}\omega + \theta_0 + h \theta_1 + \dots + h^{k-2}\theta_{k-2},
\]
where
\begin{align*}
 \omega &= (I_{-1} + qR_{-1})dt = I_{-1}dt + R_{-1}dq,\\
 \theta_0 &= R_0 dq,\ \dots,\ \theta_{k-2} = R_{k-2}dq,\\
 R_i &= \diag_{N-k+i}(0,\dots,0,\lambda_{k-2-i}).
\end{align*}
Note that $qdt=dq$ because $q=e^t$.

The following properties will be useful in the calculation of $L_+$.
\begin{prop}[\cite{guest05:_quant_d}]
 If we set $\deg h = 2$ and $\deg q = 2(N-k)$,
 the following statements hold.
 \begin{enumerate}
  \item If the $(\alpha,\beta)$-component of $L_+$ does not vanish,
	it has degree $2(\beta-\alpha)$.
  \item If the $(\alpha,\beta)$-component of $Q_i$ does not vanish,
	it has degree $2(\beta-\alpha-i)$.
  \item There is a matrix $X$ such that $Q_0 = \exp X$ and 
	the $n$-diagonal component of $X$ vanishes for $n \leq 1$.
  \item For $i \geq 1$ and $n \leq i+1$,
	the $n$-diagonal component of $Q_i$ vanishes.
 \end{enumerate}
\end{prop}
According to the above proposition, we may assume that the $Q_i$ are of
the form
\begin{align*}
 Q_0 &= I + q Q_0^1 + q^2 Q_0^2 + \dots
 = I + \sum_{\alpha \geq 1} q^\alpha Q_0^{\alpha}, \\
 Q_i &= \sum_{\alpha \geq 1} q^\alpha Q_i^{\alpha} \quad (i \geq 1),
\end{align*}
where $Q_i^\alpha$ is a constant $(i + \alpha(N-k))$-diagonal
matrix. Thus $Q_i^\alpha$ vanishes if $\alpha$ is greater than
$(N-2-i)/(N-k)$.

Before solving the equations for $L_+$, we note the following identities:
\begin{enumerate}
 \item $R_jQ_i^\alpha = 0 \quad (i,\alpha \geq 0,\ j \geq -1)$.
 \item $[Q_1,\omega] = [Q_1^1,I_{-1}]dq + \sum_{\alpha \geq 1} q^\alpha \left([Q_1^{\alpha+1},I_{-1}] + Q_1^\alpha R_{-1} \right)dq$.
\end{enumerate}

First, we consider the equation $(\L_0)$. The left hand side is
$dQ_0 = \sum_{\alpha \geq 1} \alpha q^{\alpha - 1} Q_0^{\alpha}$,
while the right hand side is 
\begin{align*}
 &Q_0(\theta_0 + [Q_1,\omega]) \\
 =\ 
 &Q_0 \Biggl( (R_0 + [Q_1^1,I_{-1}]) + \sum_{\beta \geq 1} q^{\beta}
 \Bigl( [Q_1^{\beta+1},I_{-1}] + Q_1^{\beta}R_{-1} \Bigr) \Biggr)dq \\
 =\ &(R_0 + [Q_1^1,I_{-1}])dq
 + q\Bigl(Q_0^1(R_0 + [Q_1^1,I_{-1}]) + [Q_1^2,I_{-1}]+Q_1^1R_{-1} \Bigr)dq\\
 &\hspace{69pt} + \sum_{\gamma \geq 2} q^\gamma \Bigl( [Q_1^{\gamma+1},I_{-1}]
 + Q_1^\gamma R_{-1} + Q_0^\gamma (R_0+[Q_1^1,I_{-1}]) \\
 &\hspace{133pt} + \sum_{\alpha + \beta = \gamma}
 Q_0^\alpha([Q_1^{\beta+1},I_{-1}]+Q_1^\beta R_{-1}) \Bigr)dq.
\end{align*}
Thus we have
\begin{align*}
 Q_0^1 &= R_0 + [Q_1^1,I_{-1}], \\
 2Q_0^2 &= Q_0^1(R_0 + [Q_1^1,I_{-1}]) + [Q_1^2,I_{-1}]+Q_1^1R_{-1}, \\
 \gamma Q_0^\gamma &=
 [Q_1^{\gamma},I_{-1}]
 + Q_1^{\gamma-1} R_{-1} + Q_0^{\gamma-1} (R_0+[Q_1^1,I_{-1}]) \\
 &\qquad + \sum_{\alpha,\beta \geq 1,\ \alpha + \beta = \gamma-1}
 Q_0^\alpha([Q_1^{\beta+1},I_{-1}]+Q_1^\beta R_{-1}) \quad (\gamma \geq 3).
\end{align*}

Second, we consider the equation $(\L_1)$.
The left hand side is
$dQ_1 = \sum_{\alpha \geq 1} \alpha q^{\alpha-1}Q_1^\alpha dq$,
while the terms in the right hand side are
\begin{align*}
  \theta_1 &= R_1 dq, \\
 [Q_1,\theta_0] &= \sum_{\alpha \geq 1}q^\alpha Q_1^\alpha R_0 dq \\
 [Q_2,\omega] &= [Q_2^1,I_{-1}]dq + \sum_{\alpha \geq 1} q^{\alpha}
 \bigl( [Q_2^{\alpha+1},I_{-1}] + Q_2^\alpha R_{-1} \bigr)dq, \\
 [Q_1,\omega]Q_1 &= \sum_{\alpha \geq 1}q^{\alpha}[Q_1^1,I_{-1}]Q_1^{\alpha} dq
 + \sum_{\gamma \geq 2} q^\gamma \sum_{\alpha + \beta = \gamma}
 \bigl( [Q_1^{\alpha+1},I_{-1}] + Q_1^\alpha R_{-1} \bigr) Q_1^\beta dq.
\end{align*}
Therefore
\begin{align*}
 Q_1^1 &= R_1 + [Q_2^1,I_{-1}], \\
 2Q_1^2 &= Q_1^1 R_0 + [Q_2^2,I_{-1}] + Q_2^1 R_{-1} - [Q_1^1,I_{-1}]Q_1^1, \\
 \gamma Q_1^\gamma &= Q_1^{\gamma-1}R_0 
 + [Q_2^\gamma,I_{-1}] + Q_2^{\gamma-1}R_{-1}
 - [Q_1^1,I_{-1}]Q_1^{\gamma-1} \\
 &\qquad - \sum_{\alpha,\beta \geq 1,\ \alpha + \beta = \gamma - 1}
 \bigl( [Q_1^{\alpha + 1},I_{-1}] + Q_1^{\alpha}R_{-1} \bigr)Q_1^\beta
 \quad (\gamma \geq 3).
\end{align*}

Finally, we consider the equation $(\L_i)$. The right hand side is
$dQ_i = \sum_{\alpha \geq i} \alpha q^{\alpha-1}Q_i^\alpha dq$,
and the terms in the right hand side are 
\begin{align*}
  \theta_i &= R_i dq, \\
 Q_j \theta_{i-j} &= \sum_{\alpha \geq 1} q^\alpha Q_j^\alpha R_{i-j}dq, \\
 [Q_i,\theta_0] &=  \sum_{\alpha \geq 1} q^\alpha Q_i^\alpha R_0 dq, \\
 [Q_{i+1},\omega] &= [Q_{i+1}^1,I_{-1}]dq + \sum_{\alpha \geq 1} q^\alpha
 \bigl( [Q_{i+1}^{\alpha+1},I_{-1}] + Q_{i+1}^\alpha R_{-1} \bigr)dq, \\
 [Q_1,\omega]Q_i &= \sum_{\alpha \geq 1}q^{\alpha}[Q_1^1,I_{-1}]Q_i^\alpha dq
 + \sum_{\gamma \geq 2} q^\gamma \sum_{\alpha + \beta = \gamma}
 \bigl( [Q_1^{\alpha+1},I_{-1}] + Q_1^\alpha R_{-1} \bigr) Q_i^\beta dq.
\end{align*}
Thus
\begin{align*}
 Q_i^1 &= R_i + [Q_{i+1}^1,I_{-1}], \\
 2Q_i^2 &= \sum_{j=1}^{i+1} Q_j^1 R_{i-j} + [Q_{i+1}^2,I_{-1}]
 - [Q_1^1,I_{-1}]Q_i^1, \\
 \gamma Q_i^\gamma &= 
 \sum_{j=1}^{i+1} Q_j^{\gamma-1} R_{i-j}
 +[Q_{i+1}^\gamma,I_{-1}]
 - [Q_1^1,I_{-1}]Q_i^{\gamma-1} \\
 &\qquad - \sum_{\alpha,\beta \geq 1,\ \alpha + \beta = \gamma-1}
 \bigl( [Q_1^{\alpha+1},I_{-1}] + Q_1^\alpha R_{-1} \bigr) Q_i^\beta
 \quad (\gamma \geq 3).
\end{align*}

Looking at the above identities, we see that $Q_i^\gamma$ is determined
by the following information:
\begin{enumerate}
 \item $Q_i^\alpha \quad (\alpha > \gamma)$
 \item $Q_j^\beta \quad (j < i,\ 0 \leq \beta \leq k-2)$
\end{enumerate}
Since $Q_{k-2}^1 = R_{k-2} = \diag_{N-2}(\lambda_0) = \diag_{N-2}(k!)$,
we can determine $L_+=Q_0(I + hQ_1 + \dots + h^{k-2}Q_{k-2})$ from $Q_{k-2}^1$
explicitly.

\bigskip\noindent\textbf{Example:}\quad
We apply the above results for $M_7^5$.
Its Picard-Fuchs operator is
\[
 P^{7,5} = (h\rd)^6 - 5qh^4(5\rd+4)(5\rd+3)(5\rd+2)(5\rd+1).
\]
First we calculate $\PF$:
\[
 \PF = \left(
 \begin{array}{cccccc}
  0 & 0 & 0 & 0 & 0 &   120qh^3 \\
  1/h & 0 & 0 & 0 & 0 & 1250qh^2 \\
  0 & 1/h & 0 & 0 & 0 & 4375qh \\
  0 & 0 & 1/h & 0 & 0 & 6250q \\
  0 & 0 & 0 & 1/h & 0 & 3125q/h \\
  0 & 0 & 0 & 0 & 1/h & 0 \\
 \end{array}
 \right)dt.
\]
Thus we have
\begin{align*}
 R_{-1} &= \diag_1(0,0,0,0,3125), \\
 R_{0}  &= \diag_2(0,0,0,6250), \\
 R_{1}  &= \diag_3(0,0,4375), \\
 R_{2}  &= \diag_4(0,1250), \\
 R_{3}  &= \diag_5(120).
\end{align*}

Second, we calculate $Q_i$ and $L_+$.
We can put 
\begin{align*}
 Q_0 &= I             + q Q_0^1 + q^2 Q_0^2, \\
 Q_1 &= \hspace{19.1pt} q Q_1^1 + q^2 Q_1^2, \\
 Q_2 &= \hspace{19.1pt} q Q_2^1, \\
 Q_3 &= \hspace{19.1pt} q Q_3^1.
\end{align*}
where $Q_i^\alpha$ is a constant $(i+2\alpha)$-diagonal matrix.
We can determine them in the following order:
\begin{align*}
 Q_3^1 =& R_3 = \diag_5(120), \\
 Q_2^1 =& R_2 + [Q_3^1,I_{-1}] = \diag_4(120,1130), \\
 Q_1^1 =& R_1 + [Q_2^1,I_{-1}] = \diag_3(120,1010,3245), \\
 Q_0^1 =& R_0 + [Q_1^1,I_{-1}] = \diag_2(120,890,2235,3005), \\
 Q_1^2 =& \dfrac{1}{2}\bigl( [Q_1^1,R_0] + [Q_2^1,R_{-1}]
 - [Q_1^1,I_{-1}]Q_1^1 \bigr) 
 = \diag_5(367800), \\
 Q_0^2 =& \dfrac{1}{2}\bigl( Q_0^1 \bigl(R_0+[Q_1^1,I_{-1}]\bigr)
 + [Q_1^1,R_{-1}] + [Q_1^2,I_{-1}] \bigr) 
 = \diag_4(318000,2731450).
\end{align*}
Thus we have
\begin{align*}
 Q_0
 &= I + q Q_0^1 + q^2 Q_0^2 \\
 &= I + q\diag_2(120,890,2235,3005) + q^2 \diag_4(318000,2731450) \\
 &= \left(
 \begin{array}{cccccc}
  1 & 0 & 120q & 0    & 318000q^2 & 0       \\
  0 & 1 & 0    & 890q & 0         & 2731450q^2 \\
  0 & 0 & 1    & 0    & 2235q     & 0       \\
  0 & 0 & 0    & 1    & 0         & 3005q    \\
  0 & 0 & 0    & 0    & 1         & 0       \\
  0 & 0 & 0    & 0    & 0         & 1       \\
 \end{array}
 \right),
\end{align*}
and
\begin{align*}
 L_+
 &= Q_0(I + hQ_1 + h^2Q_2 + h^3Q_3) \\
 &= \left(
 \begin{array}{cccccc}
  1 & 0 & 120q & 120qh & 120qh^2 + 318000q^2 & 120qh^3 + 757200q^2h  \\
  0 & 1 & 0    & 890q  & 1010qh              & 1130qh^2 + 2731450q^2 \\
  0 & 0 & 1    & 0     & 2235q               & 3245qh                \\
  0 & 0 & 0    & 1     & 0                   & 3005q                 \\
  0 & 0 & 0    & 0     & 1                   & 0                     \\
  0 & 0 & 0    & 0     & 0                   & 1                     \\
 \end{array}
 \right).
\end{align*}

Finally, we calculate $\hat\Omega^h$:
\begin{align*}
 \hat\Omega^h
 &= \dfrac{1}{h}Q_0 \omega (Q_0)^{-1} \\
 &= \dfrac{1}{h}Q_0 (I_{-1}+qR_{-1}) (Q_0)^{-1}dt \\
 &= \dfrac{1}{h}\left(
 \begin{array}{cccccc}
  0 & 120q & 0    & 211200q^2 & 0         & 31320000q^3 \\
  1 & 0    & 770q & 0         & 692500q^2 & 0 \\
  0 & 1    & 0    & 1345q     & 0         & 211200q^2 \\
  0 & 0    & 1    & 0         & 770q      & 0 \\
  0 & 0    & 0    & 1         & 0         & 120q \\
  0 & 0    & 0    & 0         & 1         & 0 \\
 \end{array}
 \right).
\end{align*}

We will see that
$\hat\Omega^h$ agrees with the restricted Dubrovin connection $\DUB$.

\bigskip\noindent\textbf{Example:}\quad
We apply the above results for $M_5^4$.
\[
 P^{5,4} = (h\rd)^4 - 4qh^3(4\rd+3)(4\rd+2)(4\rd+1).
\]
First we calculate $\PF$:
\[
 \PF = \left(
 \begin{array}{cccccc}
   0 & 0  & 0 & 24qh^2 \\
  1/h & 0 & 0 & 176qh \\
  0 & 1/h & 0 & 384q \\
  0 & 0 & 1/h & 256q/h \\
 \end{array}
 \right)dt.
\]
Thus we have
\begin{align*}
 R_{-1} &= \diag_0(0,0,0,256), \\
 R_{0}  &= \diag_1(0,0,384), \\
 R_{1}  &= \diag_2(0,176), \\
 R_{2}  &= \diag_3(24). \\
\end{align*}

Second, we calculate $Q_i$ and $L_+$.
We can put 
\begin{align*}
 Q_0 &= I             + q Q_0^1 + q^2 Q_0^2 + q^3 Q_0^3, \\
 Q_1 &= \hspace{19.1pt} q Q_1^1 + q^2 Q_1^2, \\
 Q_2 &= \hspace{19.1pt} q Q_2^1.
\end{align*}
where $Q_i^\alpha$ is a constant $(i+\alpha)$-diagonal matrix.
We can determine them in the following order:
\begin{align*}
 Q_2^1 &= R_2 = \diag_3(24), \\
 Q_1^1 &= R_1 + [Q_2^1,I_{-1}] = \diag_2(24,152), \\
 Q_0^1 &= R_0 + [Q_1^1,I_{-1}] = \diag_1(24,128,232), \\
 Q_1^2 &= \dfrac{1}{2}\bigl( [Q_1^1,R_0] + [Q_2^1,R_{-1}]
 - [Q_1^1,I_{-1}]Q_1^1 \bigr) 
 = \diag_3(5856), \\
 Q_0^2 &= \dfrac{1}{2}\bigl( Q_0^1 \bigl(R_0+[Q_1^1,I_{-1}]\bigr)
 + [Q_1^2,I_{-1}] + Q_1^1R_{-1} \bigr) 
 = \diag_2(4464,31376), \\
 Q_0^3 &=
 \dfrac{1}{3}\left([Q_1^3,I_{-1}] + Q_1^2R_{-1}
 + Q_0^2(R_0+[Q_1^1,I_{-1}])
 + Q_0^1([Q_1^2,I_{-1}]+Q_1^1R_{-1})
 \right) \\
 &= \diag_3(1109376).
\end{align*}

Thus we have
\begin{align*}
 Q_0
 &= I + q Q_0^1 + q^2 Q_0^2 + q^3 Q_0^3 \\
 &= I + q\diag_1(24,128,232) + q^2\diag_2(4464,31376) + q^3\diag_3(1109376)
 \\
 &= \left(
 \begin{array}{cccccc}
  1 & 24q & 4464q^2 & 1109376q^3 \\
  0 & 1   & 128q    & 31376q^2   \\
  0 & 0   & 1       & 232q       \\
  0 & 0   & 0       & 1          \\
 \end{array}
 \right),
\end{align*}
and
\begin{align*}
 L_+
 &= Q_0(I + hQ_1 + h^2Q_2) \\
 &= \left(
 \begin{array}{cccccc}
  1 & 24q & 24qh + 4464q^2 & 24qh^2 + 9504q^2h + 1109376q^3 \\
  0 & 1   & 128q           & 152qh + 31376q^2               \\
  0 & 0   & 1              & 232q                           \\
  0 & 0   & 0              & 1                              \\
 \end{array}
 \right).
\end{align*}

Finally, we calculate $\hat\Omega^h$:
\[
 \hat\Omega^h
 = \dfrac{1}{h}Q_0 \omega (Q_0)^{-1} \\
 = \dfrac{1}{h}\left(
 \begin{array}{cccccc}
  24q & 3888q^2 & 504576q^3 & 18323712q^4 \\
  1   & 104q    & 13600q^2  & 504576q^3   \\
  0   & 1       & 104q      & 3888q^2     \\
  0   & 0       & 1         & 24q         \\
 \end{array}
 \right).
\]
In fact, $\hat\Omega^h$ does not agree with the Dubrovin connection $\DUB$.
But we will see that the modified connection
\[
 -\dfrac{4!q}{h}Idt + \hat\Omega^h
 = \dfrac{1}{h}\left(
 \begin{array}{cccccc}
  0 & 3888q^2 & 504576q^3 & 18323712q^4 \\
  1 & 80q     & 13600q^2  & 504576q^3   \\
  0 & 1       & 80q       & 3888q^2     \\
  0 & 0       & 1         & 0           \\
 \end{array}
 \right)
\]
agrees with $\DUB$.

The above algorithm can easily be implemented in Maple
\footnote{A Maple program can be found at \textsf{http://blueskyproject.net/sakai/via\_dmod/} .},
Mathematica etc. It is more elementary than the method of
\cite{jinzenji02:_gauss_manin}, \cite{bertram05:_new_gromov_witten}.

\section{Adapted gauge group}

There are three important ingredients of the theory, i.e.,
$\D$-modules (with adapted basis), ``adapted'' systems of
differential equations, and ``adapted'' (flat) connections.
These are closely related to each other.

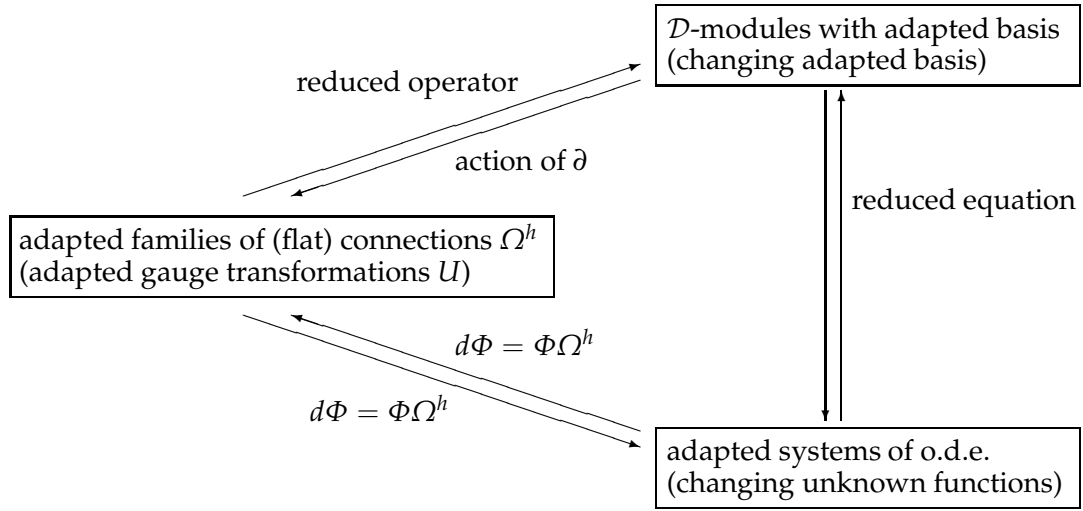
\begin{center}
\begin{picture}(430,200)(,)
 \put(2,82){\framebox(205,30)}
 \put(5,100){adapted families of (flat) connections $\Omega^h$}
 \put(5,88){(adapted gauge transformations $U$)}
 \put(90,120){\vector(3,1){150}}
 \put(110,160){reduced operator}
 \put(240,164){\vector(-3,-1){132}}
 \put(170,130){action of $\rd$}
 \put(246,162){\framebox(160,30)}
 \put(250,180){$\D$-modules with adapted basis}
 \put(250,168){(changing adapted basis)}
 \put(310,160){\vector(0,-1){125}}
 \put(320,115){reduced equation}
 \put(316,35){\vector(0,1){125}}
 \put(246,2){\framebox(160,30)}
 \put(250,20){adapted systems of o.d.e.}
 \put(250,8){(changing unknown functions)}
 \put(90,75){\vector(3,-1){150}}
 \put(170,60){$d\Phi = \Phi\Omega^h$}
 \put(240,31){\vector(-3,1){132}}
 \put(115,35){$d\Phi = \Phi\Omega^h$}
\end{picture}
\end{center}

Let $t$ be a coordinate function on $\C$ and $q=e^t$.
We consider families of connections and the gauge group on the trivial
bundle $\C \times \C^{N-1} \rightarrow \C$.
The space of connection $1$-forms on the bundle is the space of
$\End(\C^{N-1})$-valued functions on $\C$ and the gauge group is 
the space of $GL(\C^{N-1})$-valued functions $\C$.
Therefore if we think of $h$ as the loop parameter,
then the space of families of gauge transformations is identified
with $\Lambda GL(\C^{N-1})$-valued functions on $\C$.

A matrix-valued function $A$ in $q$ and $h$ is called homogeneous if
the $n$-diagonal component of $A$ has degree $2n$ for each $n$.
A $\Lambda GL(\C^{N-1})$-valued function $U$ on $\C$ is called 
adapted if $U$ satisfies following properties.
\begin{description}
 \item[(P)] $U$ is a $\Lambda GL(\C^{N-1})$-valued polynomial in $q$ and
	    $h$.
 \item[(H)] $U$ is homogeneous.
 \item[(I)] $U|_{q=0} = I$. \textsf{(initial condition)}
\end{description}
Note that the inverse of an adapted $\Lambda GL(\C^{N-1})$-valued
function is also adapted. Thus the space $\agg$ of adapted
$\Lambda GL(\C^{N-1})$-valued functions is a subgroup of the gauge group
which is called the adapted gauge group. The adapted gauge group
acts on the space of families of connection $1$-forms as same as the
gauge group. 
\[
 U^*\Omega^h = U^{-1}dU + U^{-1}\Omega^hU,
\]
where $U \in \agg$ and $\Omega^h$ is a family of connection $1$-forms.

\subsection{An adapted family of connection $1$-forms}

A family of connection $1$-forms $\Omega^h$ is called adapted
if $\Omega^h$ satisfies following properties.
\begin{description}
 \item[(P)] There is an $\End(\C^{N-1})$-valued polynomial
	    $R=R^h(q)$ in $q$ and $h$ such that
	    $\Omega^h = \dfrac{1}{h}R^h(q)dt$.
 \item[(H)] $\dfrac{1}{h}R^h(q)$ is homogeneous.
 \item[(I)] $R^h(0) = I_{-1}$.
 \item[(N)] The $(-1)$-diagonal component of $R^h(q)$ is $I_{-1}$.
	  \textsf{(normalization)}
\end{description}

We denote by $\afc$ the space of adapted families of connection $1$-forms.
Since $\rd = \rd/\rd t = q \rd/\rd q$ preserves degree,
for $U \in \agg$ and $\Omega^h \in \afc$ the family of connection $1$-forms
$U^*\Omega^h$ is also adapted. In other words, the adapted gauge
group $\agg$ acts on the space $\afc$ of adapted families of
connection $1$-forms.

\begin{thm}[Uniqueness]
 \label{thm:uniqueness}
 If $\Omega^h_1, \Omega^h_2 \in \afc$ are $\dfrac{1}{h}$-linear
 (i.e. $h\Omega^h_1, h\Omega^h_2$ are independent of $h$)
 and adapted gauge equivalent, then $\Omega^h_1 = \Omega^h_2$.
\end{thm}

\proof
Let $U \in \agg$ such that $\Omega^h_2 = U^*\Omega^h_1$.
Since $U$ is adapted, $U$ can be written as
\[
 U = \tilde{Q}_0
 \left( I + h\tilde{Q}_1 + h^2\tilde{Q}_1 + \dots + h^p\tilde{Q}_p \right)
\]
for some integer $p$. Here $\tilde{Q}_i$ is independent of $h$.
Note that $U|_{q=0} = I$ implies $\tilde{Q}_0|_{q=0} = I$ and 
$\tilde{Q}_i|_{q=0} = 0$. Comparing coefficients of powers of $h$ in the
identity
\[
 dU = U\Omega^h_2 - \Omega^h_1U
\]
gives the system:
\begin{align*}
 \Omega^h_1 &= \tilde{Q}_0 \Omega^h_2 \tilde{Q}_0^{-1}, \\
 d\tilde{Q}_0 &= \tilde{Q}_0[\tilde{Q}_1,h\Omega^h_2], \\
 d\tilde{Q}_\alpha &= [\tilde{Q}_{\alpha+1},h\Omega^h_2]
 - [\tilde{Q}_1,h\Omega^h_2]\tilde{Q}_\alpha \qquad (\alpha = 1,\cdots,p).
\end{align*}
Here $\tilde{Q}_{p+1} = 0$.
By the following lemma \ref{lem:hom_sol},
$\tilde{Q}_{\alpha+1}=0$ implies that $\tilde{Q}_\alpha$ vanishes for
$\alpha=1, \dots, p$ since $[\tilde{Q}_1,h\Omega^h_2]$ is homogeneous.
$\tilde{Q}_1=0$ implies that $d\tilde{Q}_0$ also vanishes.
Since $\tilde{Q}_0|_{q=0} = I$, we conclude that $\tilde{Q}_0 = I$.
Thus $U = I$.
\qed

\begin{lem}
 \label{lem:hom_sol}
 Let $A=A(q)$ be a homogeneous $\End(\C^{N-1})$-valued polynomial in $q$
 such that $A(0) = 0$. If an $\End(\C^{N-1})$-valued polynomial
 $X = X(q)$ satisfies the following two conditions, then $X=0$.
 \begin{enumerate}
  \item $dX = AX$.
  \item There is an integer $m$ such that the $n$-diagonal component of
	$X$ has degree $2n+2m$ for each $n$.
 \end{enumerate}
\end{lem}

\proof
Let us write $A$ and $X$ as finite sums as follows:
\[
 A = \sum_{\alpha \geq 0} q^\alpha A_\alpha, \quad
 X = \sum_{\alpha \geq 0} q^\alpha X_\alpha.
\]
Note that the initial condition on $A$ implies $A_0 = 0$.
The equation can be written as
\[
 \gamma X_\gamma = \sum_{\alpha + \beta = \gamma}A_\alpha X_\beta
 \qquad(\gamma = 0,1,2,\dots).
\]
We have $X_0=0$ since $0 = X_0 + A_0$.
If $X_\beta=0$ for $\beta = 0, \dots, \gamma$, then the identity
\[
 (\gamma+1)X_{\gamma+1}
 = \sum_{\beta = 0}^\gamma A_{(\gamma+1)-\beta} X_\beta + A_0 X_\gamma
\]
implies $X_{\gamma+1}=0$. 
By induction on $\gamma$, we conclude that $X_\gamma=0$ for all $\gamma$,
i.e., $X=0$.
\qed

\subsection{An adapted system of ordinary differential equations}
\label{subsection:adapted_system}

We consider the system $d\Phi = \Phi\Omega^h$ (with parameter $h$) of 
ordinary differential equations,
where $\Phi = (\varphi_0,\dots,\varphi_{N-2})$.
The system $d\Phi = \Phi\Omega^h$ is called adapted if
$\Omega^h$ is an adapted family of connection $1$-forms.

If we introduce new unknown functions $\Psi = \Phi U$ with $U \in \agg$,
then the system of o.d.e. is equivalent to a new adapted system of
o.d.e. $d\Psi = \Psi (U^*\Omega^h)$.

An adapted system can be reduced to an ordinary differential equation.
There is an $\End(\C^{N-1})$-valued polynomial 
$R=R^h(q)=\left(r_{\alpha,\beta}\right)_{0 \leq \alpha,\beta \leq N-2}$
in $q$ and $h$ such that $\Omega^h = \dfrac{1}{h}R^h(q)dt$.
Since $r_{\alpha,\beta} = 0$ for $\alpha > \beta+1$,
the system is written as
\[
 h\dfrac{\rd \varphi_\beta}{\rd t}
 = \sum_{\alpha=0}^{\beta+1} r_{\alpha,\beta} \varphi_{\alpha}
 \quad (\beta = 0,\dots,N-2),
\]
where $\varphi_{N-1} = 0$. Because $r_{\beta+1,\beta}=1$
($\beta = 0,\dots,N-3$),
\[
 \varphi_{\beta+1}
 = h\dfrac{\rd \varphi_\beta}{\rd t}
 - \sum_{\alpha=0}^\beta r_{\alpha,\beta}\varphi_{\alpha}
 \quad (\beta = 0,\dots,N-2).
\]
For each $\beta$, $\varphi_\beta$ can be written in terms of $\varphi_0$
and its derivatives:
\begin{description}
 \item[(P)] There exist polynomials $\sigma_{\beta,\gamma}$ in $q$ and
	    $h$ so that
	    $\varphi_\beta = \displaystyle\sum_{\gamma=0}^\beta \sigma_{\beta,\gamma} h^\gamma\dfrac{\rd^\gamma \varphi_0}{\rd t^\gamma}$.
 \item[(H)] For each $\beta$ and $\gamma$,
	    $\sigma_{\beta,\gamma}$ is a homogeneous polynomial
	    of degree $2(\beta-\gamma)$.
 \item[(I)] If $\beta > \gamma$, then $\sigma_{\beta,\gamma}|_{q=0}=0$.
 \item[(N)] $\sigma_{\beta,\beta}=1$.
\end{description}

If the above conditions hold for $\varphi_0, \dots, \varphi_\beta$, then
\begin{align*}
 \varphi_{\beta+1}
 &= h\dfrac{\rd}{\rd t} \left(\sum_{\gamma=0}^\beta \sigma_{\beta,\gamma}h^\gamma \dfrac{\rd^\gamma\varphi_0}{\rd t^\gamma}\right)
 - \sum_{\alpha=0}^\beta r_{\alpha,\beta}
 \sum_{\gamma=0}^\alpha \sigma_{\alpha,\gamma}h^\gamma \dfrac{\rd^\gamma\varphi_0}{\rd t^\gamma} \\
 &= \sigma_{\beta,\beta}h^{\beta+1} \dfrac{\rd^{\beta+1} \varphi_0}{\rd t^{\beta+1}}
 + \sum_{\gamma=1}^{\beta}
 \left(\sigma_{\beta,\gamma-1} + qh\dfrac{\rd \sigma_{\beta,\gamma}}{\rd q}
 - \sum_{\alpha=\gamma}^{\beta} r_{\alpha,\beta} \sigma_{\alpha,\gamma}
 \right) h^\gamma \dfrac{\rd^\gamma \varphi_0}{\rd t^\gamma} \\
 &\hspace{175pt}
 + \left( qh\dfrac{\rd \sigma_{\beta,0}}{\rd q}
 - \sum_{\alpha=0}^{\beta} r_{\alpha,\beta}\sigma_{\alpha,0} \right) \varphi_0.
\end{align*}
Since $\deg r_{\alpha,\beta} = 2(\beta-\alpha+1)$ and
$\deg \sigma_{\beta,\gamma} = 2(\beta-\alpha+1)$,
the conditions hold for $\varphi_{\beta+1}$.

In particular $\varphi_{N-1} = 0$. Therefore $\varphi_0$ satisfies the
following o.d.e.:
\[
 \left( (h\rd)^{N-1} + \sigma_{N-1,N-2}(h\rd)^{N-2} + \dots +
 \sigma_{N-1,1}(h\rd) + \sigma_{N-1,0} \right)\varphi_0 = 0.
\]
The above o.d.e is called the reduced equation of the system.

\subsection{A $\D$-module and an adapted basis}

We consider a $\D$-module $\M^h = \D/(P)$ of rank $N-1$ for some
differential operator $P \in \D$.
We assume that there exists homogeneous polynomials $a_\alpha$ of degree
$2\alpha$ such that
$P = (h\rd)^{N-1} + \sum_{\alpha=1}^{N-1} a_\alpha(h\rd)^{N-1-\alpha}$.
Let $P_0,\dots,P_{N-2} \in \D$ be differential operators such that
$[P_0],\dots,[P_{N-2}]$ form a $\Lambda(h)$-basis of a $\M^h$. We say
that $[P_0],\dots,[P_{N-2}]$ form an adapted ($\Lambda(h)$-)basis if
the differential operators satisfy the following properties. 
\begin{description}
 \item[(P)] There exist polynomials $c_{\beta,\gamma}$ in $q$ and $h$ so
	    that
	    $P_\beta = \sum_{\gamma=0}^\beta c_{\beta,\gamma}(h\rd)^\gamma$.
 \item[(H)] For each $\beta$ and $\gamma$,
	    $c_{\beta,\gamma}$ is a homogeneous polynomial
	    of degree $2(\beta-\gamma)$.
 \item[(I)] If $\beta > \gamma$, then $c_{\beta,\gamma}|_{q=0}=0$.
 \item[(N)] $c_{\beta,\beta}=1$.
\end{description}

A $\D$-module $\M^h = \D/(P)$ with an adapted basis
$[P_0],\dots,[P_{N-2}]$ defines a family of (flat) connection
$1$-forms $\Omega^h$ as in section \ref{section:Birkhoff}.
\[
 \Omega^h = \dfrac{1}{h}R^h(q)dt
 \quad\mbox{where}\quad
 h\rd ([P_0],\dots,[P_{N-2}]) = ([P_0],\dots,[P_{N-2}]) R^h(q).
\]
Let $[P'_0],\dots,[P'_{N-2}]$ be another adapted basis of $\M^h$.
The adapted conditions imply that there exists an adapted gauge
transformation $U \in \agg$ such that 
\[
 (P'_0, \dots, P'_{N-2}) = (P_0, \dots, P_{N-2})U.
\]
The family $\Omega'^h$ of connection $1$-forms associated with the
adapted basis $[P'_0],\dots,[P'_{N-2}]$ agrees with $U^*\Omega^h$.
In fact,
\begin{align*}
 h\rd ([P'_0],\dots,[P'_{N-2}])
 &= h\rd \bigl(([P_0],\dots,[P_{N-2}])\bigr)U
 + ([P_0],\dots,[P_{N-2}])h\rd U \\
 &= ([P_0],\dots,[P_{N-2}]) R^h(q)U
 + ([P_0],\dots,[P_{N-2}])h\rd U \\
 &= ([P'_0],\dots,[P'_{N-2}]) \left( U^{-1}R^h(q)U + U^{-1}h\rd U \right).
\end{align*}
Therefore
\[
 \Omega'^h = \dfrac{1}{h}\left( U^{-1}R^h(q)U + U^{-1}h\rd U \right)dt
 = U^*\Omega^h.
\]

Conversely, an adapted family $\Omega^h$ of connection $1$-forms
defines a $\D$-module $\M^h$ and an adapted
$\Lambda(h)$-basis as follows. There is an $\End(\C^{N-1})$-valued
polynomial 
$R=R^h(q)=\left(r_{\alpha,\beta}\right)_{0 \leq \alpha,\beta \leq N-2}$
in $q$ and $h$ such that $\Omega^h = \dfrac{1}{h}R^h(q)dt$.
The differential operators $P_0, \dots, P_{N-1}$ are defined inductively
as follows.
\[
 P_0 = 1, \qquad
 P_{\beta+1} = (h\rd)P_{\beta}
 - \sum_{\alpha = 0}^{\beta} r_{\alpha,\beta} P_{\alpha}
 \quad (\beta = 0, \cdots, N-2).
\]
We define a $\D$-module $\M^h$ by $\D/(P_{N-1})$
and we call the operator $P_{N-1}$ the reduced operator of $\Omega^h$.
As in subsection \ref{subsection:adapted_system},
$[P_0],\dots,[P_{N-2}]$ form an adapted $\Lambda(h)$-basis of $\M^h$.
By the definition of the operators $P_0, \dots, P_{N-2}$, we have
$h\rd ([P_0],\dots,[P_{N-2}]) = ([P_0],\dots,[P_{N-2}]) R^h(q)$.

\begin{lem}
 \label{lem:same_reduced}
 Let $\M^h$ be a $\D$-module $\D/(P)$ of rank $N-1$ and
 $\Omega^h$ an adapted family associated with an adapted basis
 $[P_0],\dots,[P_{N-2}]$. Then the reduced operator of $\Omega^h$ agrees
 with $P$. In particular, the reduced operator is independent of the
 choice of the adapted basis.
\end{lem}

\proof
Let $P'=(h\rd)P_{N-2}-\sum_{\alpha=0}^{N-2}r_{\alpha,N-2}(h\rd)^\alpha$
be the reduced operator of $\Omega^h$. By the definition of
$\Omega^h = \dfrac{1}{h}\left(r_{\alpha,\beta}\right)_{0 \leq \alpha,\beta \leq N-2}dt$, we have
\[
 \left[ (h\rd)P_{N-2} \right] = 
 \left[ \sum_{\alpha=0}^{N-2} r_{\alpha,N-2} (h\rd)^\alpha \right],
 \quad\mbox{i.e.}\quad [P'] = 0.
\]
Since the monic polynomial $P' \in \Lambda(h)[h\rd]$ has order $N-1$,
$P'$ must agree with $P$.
\qed

\begin{thm}[Equivalence]
 \label{thm:equivalence}
 Let $\Omega^h_1, \Omega^h_2 \in \afc$. 
 Then $\Omega^h_1$ and $\Omega^h_2$ have a same reduced operator
 if and only if
 $\Omega^h_1$ and $\Omega^h_2$ are adapted gauge equivalent.
\end{thm}

\proof
First, we assume that $\Omega^h_1$ and $\Omega^h_2$ have a
same reduced operator. Let $P_0^\alpha, \dots, P_{N-1}^\alpha$ be the
operators associated with $\Omega^h_\alpha$ for $\alpha = 1,2$. By the
assumption, $P = P_{N-1}^1 = P_{N-1}^2$. Since two basis
$[P_0^1], \dots, [P_{N-2}^1]$ and $[P_0^2], \dots, [P_{N-2}^2]$
are adapted basis of the same $\D$-module $\M^h = \D/(P)$,
there exists an adapted gauge transformation $U \in \agg$
such that $([P_0^2],\dots,[P_{N-2}^2]) = ([P_0^1],\dots,[P_{N-2}^1])U$.
If $R^h_\alpha(q)$ is a matrix-valued polynomial in $q$ and $h$ such
that $\Omega^h_\alpha = \dfrac{1}{h}R^h_\alpha(q)$, then
$h\rd ([P_0^\alpha], \dots, [P_{N-2}^\alpha]) = ([P_0^\alpha], \dots, [P_{N-2}^\alpha])R^h_\alpha(q)$ for $\alpha = 1,2$.
Therefore we conclude that $\Omega^h_2 = U^*\Omega^h_1$.

Next, we assume that there exists an adapted gauge
transformation $U \in \agg$ such that $\Omega^h_1$ and $\Omega^h_2$.
Let $\M^h$ with $[P_0^1], \dots, [P_{N-2}^1]$ be the $\D$-module with
the adapted basis defined by $\Omega^h_1$ and $P_{N-1}^1$ the reduced
operator of $\Omega^h_1$. Define an adapted basis
$[P_0^2], \dots, [P_{N-2}^2]$ by
$(P_0^2, \dots, P_{N-2}^2) = (P_0^1, \dots, P_{N-2}^1)U$.
Then the adapted family $\Omega^h_2$ agrees with the adapted
family defined by the adapted basis $[P_0^2], \dots, [P_{N-2}^2]$ of
$\M^h = \D/(P_{N-1}^1)$. According to the lemma \ref{lem:same_reduced},
the reduced operator of $\Omega^h_2$ agrees with $P_{N-1}^1$.
\qed

\section{Relation between Birkhoff factorization and Jinzenji's results}
\label{section:Jinzenji}

Recall the quantum differential system (with parameter $h$)
for $M_N^k$:
\begin{align*}
 &h\dfrac{\rd \psi_{N-2-m}}{\rd t} = \psi_{N-1-m}(t)
 + \sum_{d \geq 1} L_m^d q^d \psi_{N-1-m-d(N-k)}(t),\\
 &h\dfrac{\partial \psi_{N-2}}{\partial t} = 
 \sum_{d \geq 1} L_0^d q^d \psi_{N-1-d(N-k)}(t),
\end{align*}
(where $m = 1, \dots, N-2$).  Note that $\deg h = 2$ and
$\deg q = 2(N-k)$. We can write the above system of o.d.e. with the
restricted Dubrovin connection $\DUB \in \afc$:
\[
 d\Psi = \Psi \DUB,
 \quad\mbox{where}\quad
 \Psi = (\psi_0,\dots,\psi_{N-2}).
\]
Theorem \ref{thm:Givental} says the reduced operator of $\DUB$
at $h=1$ agrees with the Picard-Fuchs operator
$\rd^{N-1} - kq(k\rd+(k-1))\cdots(k\rd+1)$ if $N-k \geq 2$.
We shall now give the proof of proposition \ref{prop:reduction},
which states that the reduced operator of $\DUB$ agrees with the 
the following operator:
\[
 (h\rd)^{N-1} - kqh^{k-1} \bigl(k\rd + (k-1)\bigr)\dots\bigl(k\rd+1\bigr).
\]

\proofof{proposition \ref{prop:reduction}}
Let $P = (h\rd)^{N-1} + c_{N-1,N-2}(h\rd)^{N-2} + \dots + c_{N-1,0}$ be
the reduced operator of $\DUB$. Since $P|_{h=1}$ agrees with the
Picard-Fuchs operator, we have
\[
 c_{N-1,\gamma}|_{h=1} =
 \begin{cases}
  -\lambda_\gamma q & (0 \leq \gamma \leq k-1), \\
  0 & (k \leq \gamma \leq N-2).
 \end{cases}
\]
Since $c_{N-1,\gamma}$ is a homogeneous polynomial of degree
$2(N-1-\gamma)$, we have
\[
 c_{N-1,\gamma} =
 \begin{cases}
  -\lambda_\gamma q h^{k-1-\gamma}& (0 \leq \gamma \leq k-1), \\
  0 & (k \leq \gamma \leq N-2).
 \end{cases}
\]
Therefore
\begin{align*}
 P 
 &= (h\rd)^{N-1} - \lambda_{k-1}qh^{k-1}\rd^{k-1}
 - \lambda_{k-2}qh^{k-1}\rd^{k-2} - \dots - \lambda_1qh^{k-1}\rd
 - \lambda_0 qh^{k-1} \\
 &= (h\rd)^{N-1} - kqh^{k-1}
 \bigl(k\rd + (k-1)\bigr)\dots\bigl(k\rd+1\bigr).
\end{align*}
\qed

\begin{thm}
 Let $\Omega^h \in \afc$ be an adapted family of connection $1$-forms
 whose reduced operator agrees with
 $P^{N,k}=(h\rd)^{N-1}-kqh^{k-1} \bigl(k\rd+(k-1)\bigr)\dots\bigl(k\rd+1\bigr)$.
 If $\Omega^h$ is $h^{-1}$-linear,  then $\Omega^h = \hat\Omega^h$.
\end{thm}
Note that we have the adapted family $\PF \in \afc$,
whose reduced operator agrees with $P^{N,k}$.
Using a matrix-valued function $L$ which satisfies $\PF=L^{-1}dL$,
$\hat\Omega^h$ is defined as $(L_-)^{-1}dL_-$ in section
\ref{section:Birkhoff}. Here $L_-$ is the first factor of the Birkhoff
factorization of $L=L_-L_+$. Since $L_+ \in \agg$, $\PF$ and
$\hat\Omega^h$ are adapted gauge equivalent.

\proof
Note that $\hat\Omega^h$ is adapted and $h^{-1}$-linear.
According to theorem \ref{thm:equivalence},
$\PF$ and $\hat\Omega^h$ has the same reduced operator $P^{N,k}$,
because $\PF$ and $\hat\Omega^h$ are adapted gauge equivalent.
Therefore $\hat\Omega^h$ satisfies the conditions of the theorem.

If $\Omega^h$ satisfies the conditions of the theorem,
then theorem \ref{thm:equivalence} says $\Omega^h$ and $\hat\Omega^h$ are
adapted gauge equivalent. Moreover the fact that 
$\Omega^h$ and $\hat\Omega^h$ are $h^{-1}$-linear implies 
$\Omega^h = \hat\Omega^h$ because of theorem \ref{thm:uniqueness}.
\qed

However the reduced operator of the quantum differential system for
$M_N^k$ differs for the cases $N-k \geq 2$ and $N-k=1$,
Jinzenji considered an adapted and $h^{-1}$-linear family
$\JIN \in \afc$ whose reduced operator agrees with $P^{N,k}$
and he named the coefficients of $\JIN$ the virtual structural constants.
Moreover he gave explicit formula for $\JIN$ in
\cite{jinzenji02:_gauss_manin}.
Jinzenji's explicit formulae guarantee the existence of $\JIN$.
Since the adapted family $\JIN$ automatically satisfies the
conditions of the above theorem, $\JIN$ agrees with $\hat\Omega^h$.

\begin{cor}
 The adapted family $\JIN$ agrees with $\hat\Omega^h$.
\end{cor}

In the case $N-k \geq 2$, $\JIN$ a priori agrees with the restricted
Dubrovin connection $1$-form $\DUB$, hence so does $\hat\Omega^h$.

In the case $N-k=1$, the Dubrovin connection $1$-form has different
reduced operator from the case of $N-k \geq 2$.

\begin{thm}[Givental]
 Let $S = \exp\left(-\frac{(N-1)!q}{h}\right)I$. The quantum
 differential system for $M_N^{N-1}$ can be written as
 \[
  d\Psi = \Psi (S^*\JIN).
 \]
\end{thm}
Note that $S$ is not adapted. The above theorem implies 
\[
\DUB = S^*\JIN = S^*\hat\Omega^h = -\frac{(N-1)!q}{h}Idt + \hat\Omega^h. 
\]

\def\cprime{$'$}

\end{document}